\documentclass[final,onefignum,onetabnum]{siamart}
\pdfoutput=1

\usepackage{lipsum}
\usepackage{amsfonts}
\usepackage{graphicx}
\usepackage{epstopdf}
\usepackage{algorithmic}
\ifpdf
  \DeclareGraphicsExtensions{.eps,.pdf,.png,.jpg}
\else
  \DeclareGraphicsExtensions{.eps}
\fi

\newsiamremark{remark}{Remark}
\newsiamremark{hypothesis}{Hypothesis}
\crefname{hypothesis}{Hypothesis}{Hypotheses}
\newsiamthm{claim}{Claim}

\headers{Relative Schwarz Waveform Relaxation Algorithm}{Fei Wei and Anna Zhao}

\title{A Non-iterative Overlapping Schwarz Waveform Relaxation Algorithm for Wave Equation\thanks{Submitted to the editors DATE.}
}

\author{Fei Wei\thanks{NMHW Technologies Co. Ltd, Shanghai, China
  (\email{rudeway@163.com}).}
\and Anna Zhao\thanks{NMHW Technologies Co. Ltd, Shanghai, China 
  .}
}

\usepackage{amsopn}

\ifpdf
\hypersetup{
  pdftitle={A Non-iterative Overlapping Schwarz Waveform Relaxation Algorithm for Wave Equation},
  pdfauthor={Fei Wei and Anna Zhao}
}
\fi

\begin{document}

\maketitle

\begin{abstract}
The Schwarz Waveform Relaxation algorithm (SWR) exchanges the waveform of boundary value between neighbouring sub-domains, which provides a more efficient way than the other Schwarz algorithms to realize distributed computation. However, the convergence speed of the traditional SWR is slow, and various optimization strategies have been brought in to accelerate the convergence.

In this paper, we propose a non-iterative overlapping variant of SWR for wave equation, 
which is named Relative Schwarz Waveform Relaxation algorithm (RSWR). 
RSWR is inspired by the physical observation that the velocity of wave is limited, 
based on the Theory of Relativity. 
The change of value at one space point will take time span $\Delta t$ to transmit to another space point and vice versa. 
This $\Delta t$ could be utilized to design distributed numerical algorithm, as we have done in RSWR. 

During each time span, RSWR needs only 3 steps to achieve high accurate waveform, 
by using the predict-select-update strategy. 
The key for this strategy is to find the maximum time span for the waveform. 
The validation of RSWR could be proved straightfowardly.
Numerical experiments show that RSWR is accurate, and is potential to be scalable and fast.

\end{abstract}



\section{Introduction}

The Schwarz Waveform Relaxation algorithm (SWR) is a combination of the Schwarz algorithm, 
see Schwarz \cite{Schwarz1870}, and the Waveform Relaxation algorithm (WR), 
see Lelarasmee et al. \cite{Lelarasmee1982}. 
At each time step, instead of exchanging values of one time point between neighbouring subdomains, 
SWR exchanges waveforms of one time span.
SWR was introduced by Bjørhus for hyperbolic problems \cite{Bjorhus1995}.
Later SWR was studied for the heat equation by Gander and Stuart \cite{Gander1998}, 
for the wave equation by Gander and Halpern \cite{Gander2005}, 
and for the time domain maxwell equation by Courvoisier and Gander \cite{Gander2013}.
Recently, the SWR was analysed at a semi-discrete level by Al-Khaleel and Wu \cite{Khaleel2019}, 
the parareal SWR was analysed by Gander et al \cite{Gander2019}, 
and a new two level SWR was described by Gander et al \cite{Gander2021}.

We studied the distributed computation of integrated circuits since 2004, 
and found that there was a similarity between 
distributed numerical algorithm and distributed physical circuit.
The key for this similarity is the existence of the transmission line, 
which is capable to partition the physical circuit into $N$ sub-circuits. 
We had designed several numerical algorithms by migrating the mathematical model 
of the lossless transmission line into 
distributed numerical algorithms \cite{Wei2008VTM, Wei2008DTM, Wei2009WTM}. 
The effect of inserting virtual transmission line into circuit is basically equivalent to 
the Schwarz algorithm with Robin transmission condition.
However, numerical experiments showed that all these distributed iterative algorithms 
suffer from the convergence problem,
and it is difficult to choose the proper characteristic impedance for 
the virtual transmission line to accelerate the algorithms. 
 
Later, we came to realize that there might be some other way to mimic the natural distributed physical circuit, 
therefore RSWR was thought of, designed, optimized and tested. 
The main advantage of RSWR over the traditional SWR algorithms is that 
it is capable to achieve the accurate result by 3 steps, by using the predict-select-update strategy.
Therefore RSWR might be considered as a non-iterative, or direct domain decomposition algorithm.

The paper is organized as follows. The detailed algorithm of RSWR for $N=2$ is described in
\cref{sect_algorithm}, the extended algorithm for $N>2$ is explained in \cref{section_algorithm_n}, 
experimental results are shown in \cref{section_experiment}, 
and the conclusions follow in \cref{section_conclusion}.

\section{Algorithm}
\label{sect_algorithm}
This section describes RSWR for 1-dimension wave equation. 
RSWR first splits the original domain into overlapping sub-domains, 
then uses the predict-select-update strategy to calculate the solution on each sub-domain.

\subsection{Wave Equation}
The wave equation in 1-dimension is expressed as:

\begin{equation}
\label{wave_equation}
\frac{{{\partial ^2}u(x,t)}}{{\partial {t^2}}} - {a^2}\frac{{{\partial ^2}u(x,t)}}{{\partial {x^2}}} = 0
\end{equation}

where: $a > 0, x \in \Omega, \Omega  = \left[ {{X_A},{X_B}} \right], t \in [0, + \infty )$.

The initial condition for Eq. \cref{wave_equation} is:

\begin{equation}
\left\{ 
       \begin{array}{l}
             u(x,t){|_{t = 0}} = u(x,0) = 0,x \in \Omega  \\ 
             \frac{{\partial u(x,t)}}{{\partial t}}{|_{t = 0}} = 0,x \in \Omega  \\ 
      \end{array} 
\right. 
\end{equation}

The boundary condition for Eq. \cref{wave_equation} is:

\begin{equation}
\left\{ 
\begin{array}{l}
 u(x,t){|_{x={X_A}}} = {f_A}(t),t \in [0, + \infty ) \\ 
 u(x,t){|_{x={X_B}}} = {f_B}(t),t \in [0, + \infty ) \\ 
 \end{array}
\right.
\end{equation}

\begin{definition}[True Solution]
$u(x,t)$ is called the true solution for the original equation Eq. \cref{wave_equation}.
\end{definition}

\subsection{Decomposition} 
\label{subsection_split}

This subsection decomposes the original domain into 2 sub-domains 
by following the Neumann transmission condition, 
i.e. the normal derivatives (flux) at the boundary are continuous.

Split $\Omega$ into 2 overlapping sub-domains $\Omega _1$ and $\Omega _2$ as below:

\begin{equation}
\left\{
\begin{array}{l}
 {\Omega _1} = [{X_{1A}},{X_{1B}}] \\ 
 {\Omega _2} = [{X_{2A}},{X_{2B}}] \\ 
 {X_A} = {X_{1A}} < {X_{2A}} < {X_{1B}} < {X_{2B}} = {X_B} \\ 
 \end{array}
\right.
\end{equation}

The overlapping region of $\Omega _1$ and $\Omega _2$ is:
\begin{equation}
{{\rm O}_{1,2}} = {\Omega _1} \cap {\Omega _2} = [{X_{2A}},{X_{1B}}] \ne \emptyset
\end{equation}

Eq. \cref{wave_equation} is split into 2 sub-equations Eq. \cref{sub_equation_1} and Eq. \cref{sub_equation_2}:

\begin{equation}
\label{sub_equation_1}
\left\{
\begin{array}{l}
 \frac{{{\partial ^2}{p_1}(x,t)}}{{\partial {t^2}}} - {a^2}\frac{{{\partial ^2}{p_1}(x,t)}}{{\partial {x^2}}} = 0,a > 0,x \in {\Omega _1},{\Omega _1} = \left[ {{X_{1A}},{X_{1B}}} \right],t \in [0, + \infty ) \\ 
 {p_1}(x,t){|_{t = 0}} = u(x,t){|_{t = 0}} = 0,x \in {\Omega _1} \\ 
 \frac{{\partial {p_1}(x,t)}}{{\partial t}}{|_{t = 0}} = 0,x \in {\Omega _1} \\ 
 {p_1}(x,t){|_{x = {X_{1A}}}} = {f_A}(t),t \in [0, + \infty ) \\ 
 \frac{{\partial {p_1}(x,t)}}{{\partial x}}{|_{x = {X_{1B}}}} = \frac{{\partial {p_2}(x,t)}}{{\partial x}}{|_{x = {X_{1B}}}},t \in [0, + \infty ) \\ 
 \end{array}
\right.
\end{equation}

\begin{equation}
\label{sub_equation_2}
\left\{
\begin{array}{l}
 \frac{{{\partial ^2}{p_2}(x,t)}}{{\partial {t^2}}} - {a^2}\frac{{{\partial ^2}{p_2}(x,t)}}{{\partial {x^2}}} = 0,a > 0,x \in {\Omega _2},{\Omega _2} = \left[ {{X_{2A}},{X_{2B}}} \right],t \in [0, + \infty ) \\ 
 {p_2}(x,t){|_{t = 0}} = u(x,t){|_{t = 0}} = 0,x \in {\Omega _2} \\ 
 \frac{{\partial {p_2}(x,t)}}{{\partial t}}{|_{t = 0}} = 0,x \in {\Omega _2} \\ 
 \frac{{\partial {p_2}(x,t)}}{{\partial x}}{|_{x = {X_{2A}}}} = \frac{{\partial {p_1}(x,t)}}{{\partial x}}{|_{x = {X_{2A}}}},t \in [0, + \infty ) \\ 
 {p_2}(x,t){|_{x = {X_{2B}}}} = {f_B}(t),t \in [0, + \infty ) \\ 
 \end{array}
\right.
\end{equation}

\begin{definition}[Boundary Flux Input Waveform, or Boundary Input Waveform for short]

\[\frac{{\partial {p_1}(x,t)}}{{\partial x}}{|_{x = {X_{1B}}}},t \in [{T_{start}},{T_{start}} + \Delta T]\] is called the boundary input waveform for sub-domain $\Omega _1$.

\[\frac{{\partial {p_2}(x,t)}}{{\partial x}}{|_{x = {X_{2A}}}},t \in [{T_{start}},{T_{start}} + \Delta T]\] is called the boundary input waveform for sub-domain $\Omega _2$.

\end{definition}

\begin{definition}[Boundary Flux Output Waveform, or Boundary Output Waveform for short]

\[\frac{{\partial {p_1}(x,t)}}{{\partial x}}{|_{x = {X_{2A}}}},t \in [{T_{start}},{T_{start}} + \Delta T]\] is called the boundary output waveform in sub-domain $\Omega _1$;

\[\frac{{\partial {p_2}(x,t)}}{{\partial x}}{|_{x = {X_{1B}}}},t \in [{T_{start}},{T_{start}} + \Delta T]\] is called the boundary output waveform in sub-domain $\Omega _2$.

\end{definition}

\begin{definition}[Input Boundary]

$x = {X_{1B}}$ is called the input boundary for sub-domain $\Omega _1$ ;

$x = {X_{2A}}$ is called the input boundary for sub-domain $\Omega _2$ ;

\end{definition}

\begin{definition}[Output Boundary]

$x = {X_{2A}}$ is called the output boundary in sub-domain $\Omega _1$ ;

$x = {X_{1B}}$ is called the output boundary in sub-domain $\Omega _2$ ;

\end{definition}

Because overlapping region exists, for sub-domain $\Omega _i$, the input boundary and the output boundary are not the same.

\begin{definition}[Corresponding Waveform]

If the boundary input waveform of one sub-domain $\Omega _i$ has the same boundary with the boundary output waveform of its adjacent sub-domain $\Omega _j$, then these two waveforms are called corresponding waveform for each other.

\end{definition}

For example:

\[\frac{{\partial {p_1}(x,t)}}{{\partial x}}{|_{x = {X_{1B}}}},t \in [{T_{start}},{T_{start}} + \Delta T]\]

is the boundary input waveform for sub-domain $\Omega _1$, whose corresponding waveform is: 

\[\frac{{\partial {p_2}(x,t)}}{{\partial x}}{|_{x = {X_{1B}}}},t \in [{T_{start}},{T_{start}} + \Delta T]\]

which is the boundary output waveform in sub-domain $\Omega _2$ ; 

Similarly, 

\[\frac{{\partial {p_2}(x,t)}}{{\partial x}}{|_{x = {X_{2A}}}},t \in [{T_{start}},{T_{start}} + \Delta T]\]

is the boundary input waveform for sub-domain $\Omega _2$, whose corresponding waveform is

\[\frac{{\partial {p_1}(x,t)}}{{\partial x}}{|_{x = {X_{2A}}}},t \in [{T_{start}},{T_{start}} + \Delta T]\]

which is the boundary output waveform in sub-domain $\Omega _1$.

\begin{theorem}[Decomposition Theorem]
\label{split_theorem}
After decomposition of the original wave equation, 
assume that for any sub-domain ${\Omega _i},i = 1, \cdots ,n$, 
each boundary input waveform of sub-domain $\Omega _i$ is equal to 
the corresponding boundary output flux waveform in its adjacent sub-domain $\Omega _j$, 
then the solution of sub-domain is equal to the true solution of the original equation.
\end{theorem}

\cref{split_theorem} means that if the original equation is decomposed by using the Nuemann transmission condition, then the solution of the sub-domian is consistent to the solution of the original equation. 

According to \cref{split_theorem}, for Eq. \cref{sub_equation_1}, because:

\[\frac{{\partial {p_1}(x,t)}}{{\partial t}}{|_{x = {X_{1B}}}} = \frac{{\partial {p_2}(x,t)}}{{\partial x}}{|_{x = {X_{1B}}}},t \in [0, + \infty )\]

thus:

\[{p_1}(x,t) = u(x,t),x \in {\Omega _1},t \in [0, + \infty )\]

Similarly, for Eq. \cref{sub_equation_2}, because:

\[\frac{{\partial {p_2}(x,t)}}{{\partial x}}{|_{x = {X_{2A}}}} = \frac{{\partial {p_1}(x,t)}}{{\partial x}}{|_{x = {X_{2A}}}},t \in [0, + \infty )\] 

thus:

\[{p_2}(x,t) = u(x,t),x \in {\Omega _2},t \in [0, + \infty )\]


\subsection{Prediction}
\label{subsection_predict}

This subsection predicts the boundary input waveform to be zero and calculate the predictive solution for each sub-domain.

First, set the time span index $k = 1$, and set the start time:

\[T_{_{start}}^k{|_{k = 1}} = 0\]

Then, set the predictive time span:

\[\Delta \hat T_{predict}^k{|_{k = 1}} = \hat \tau \] 

$\hat \tau$ should be set as a large enough positive value.

Assume that the boundary input waveform of Eq. \cref{sub_equation_1} is 0:

\[\frac{{\partial {{\hat p}_1}(x,t)}}{{\partial x}}{|_{x = {X_{1B}}}} = 0,t \in [T_{_{start}}^k,T_{_{start}}^k + \Delta \hat T_{predict}^k],k = 1\]

and the initial condition of Eq. \cref{sub_equation_1} is equal to the true solution $u(x,t)$ of Eq. \cref{wave_equation}, then Eq. \cref{sub_equation_1} is transferred into Eq. \cref{sub_equation_predict_1}:

\begin{equation}
\label{sub_equation_predict_1}
\left\{
 \begin{array}{l}
 \frac{{{\partial ^2}{{\hat p}_1}(x,t)}}{{\partial {t^2}}} - {a^2}\frac{{{\partial ^2}{{\hat p}_1}(x,t)}}{{\partial {x^2}}} = 0,a > 0,t \in [T_{_{start}}^k,T_{_{start}}^k + \Delta \hat T_{predict}^k], k = 1\\ 
\quad \quad  x \in {\Omega _1},{\Omega _1} = \left[ {{X_{1A}},{X_{1B}}} \right] \\
{{\hat p}_1}(x,t){|_{t = T_{_{start}}^k}} = u(x,t){|_{t = T_{_{start}}^k}},x \in {\Omega _1} \\ 
 \frac{{\partial {{\hat p}_1}(x,t)}}{{\partial t}}{|_{_{t = T_{_{start}}^k}}} = \frac{{\partial u(x,t)}}{{\partial t}}{|_{_{t = T_{_{start}}^k}}},x \in {\Omega _1} \\ 
 {{\hat p}_1}({X_{1A}},t) = {f_A}(t),t \in [T_{_{start}}^k,T_{_{start}}^k + \Delta \hat T_{predict}^k] \\ 
 \frac{{\partial {{\hat p}_1}(x,t)}}{{\partial x}}{|_{x = {X_{1B}}}} = 0,t \in [T_{_{start}}^k,T_{_{start}}^k + \Delta \hat T_{predict}^k] \\ 
 \end{array}
\right.
\end{equation}

Eq. \cref{sub_equation_predict_1} is able to be solved and the solution is:

\[{\hat p_1}(x,t),x \in {\Omega _1},t \in [T_{_{start}}^k,T_{_{start}}^k + \Delta \hat T_{predict}^k],k = 1\]

Similarly, assume that the boundary input waveform of Eq. \cref{sub_equation_2} is 0 and the initial condition is equal to the true solution $u(x,t)$ of Eq. \cref{wave_equation}, then Eq. \cref{sub_equation_2} is transferred into Eq. \cref{sub_equation_predict_2}:

\begin{equation}
\label{sub_equation_predict_2}
\left\{
\begin{array}{l}
 \frac{{{\partial ^2}{{\hat p}_2}(x,t)}}{{\partial {t^2}}} - {a^2}\frac{{{\partial ^2}{{\hat p}_2}(x,t)}}{{\partial {x^2}}} = 0, a > 0, t \in [T_{_{start}}^k,T_{_{start}}^k + \Delta \hat T_{predict}^k], k = 1\\
\quad \quad x \in {\Omega _2},{\Omega _2} = \left[ {{X_{2A}},{X_{2B}}} \right]\\
 {{\hat p}_2}(x,t){|_{t = T_{_{start}}^k}} = u(x,t){|_{t = T_{_{start}}^k}},x \in {\Omega _2} \\ 
 \frac{{\partial {{\hat p}_2}(x,t)}}{{\partial t}}{|_{_{t = T_{_{start}}^k}}} = \frac{{\partial u(x,t)}}{{\partial t}}{|_{_{t = T_{_{start}}^k}}},x \in {\Omega _2} \\ 
 \frac{{\partial {{\hat p}_2}(x,t)}}{{\partial x}}{|_{x = {X_{2A}}}} = 0,t \in [T_{_{start}}^k,T_{_{start}}^k + \Delta \hat T_{predict}^k] \\ 
 {{\hat p}_2}(x,t){|_{x = {X_{2B}}}} = {f_B}(t),t \in [T_{_{start}}^k,T_{_{start}}^k + \Delta \hat T_{predict}^k] \\ 
 \end{array}
\right.
\end{equation}

Eq. \cref{sub_equation_predict_2} is able to be solved and the solution is:

\[{\hat p_2}(x,t),x \in {\Omega _2},t \in [T_{_{start}}^k,T_{_{start}}^k + \Delta \hat T_{predict}^k],k = 1\]

\begin{definition}[Predictive Solution]

Assume that each boundary input waveform of sub-domain $\Omega _i$ is 0, therefore the solution ${\hat p_i}(x,t)$ of the sub-equation for $\Omega _i$ is called the predictive solution of $\Omega _i$.

\end{definition}

As the result,  

\[{\hat p_1}(x,t), x \in \Omega _1, t \in [T_{_{start}}^k,T_{_{start}}^k + \Delta \hat T_{predict}^k]\] 

is called the predictive solution of sub-domain $\Omega _1$;

\[{\hat p_2}(x,t), x \in \Omega _2, t \in [T_{_{start}}^k,T_{_{start}}^k + \Delta \hat T_{predict}^k]\] 

is called the predictive solution of sub-domain $\Omega _2$.


\subsection{Selection}
\label{subsection_select}

This subsection selects the maximum waveform time span by comparing the predictive solutions of adjacent sub-domains on the overlapping region.

\begin{definition}[Maximum Waveform Time Span]

In the overlapping region ${{\rm O}_{1,2}} = {\Omega _1} \cap {\Omega _2}$, define the maximum point time span \newline $\Delta T_{\max }^k(x)$ as:

For $x \in {\Omega _1} \cap {\Omega _2}$, $\Delta T_{\max }^k(x) = \max (\Delta {T^k}(x))$, where $\Delta {T^k}(x)$ satisfies $\forall t \in [T_{start}^k,T_{start}^k + \Delta {T^k}(x)]$, ${{\hat p}_1}(x,t) = {{\hat p}_2}(x,t)$, $k=1$.

Then define maximum waveform time span $\Delta T_{\max }^k$ as:

$\Delta T_{\max }^k = \max (\Delta T_{\max }^k(x)),\forall x \in {\Omega _1} \cap {\Omega _2},k = 1$

\end{definition}

\begin{theorem}[Prediction Validation Theorem]
\label{predict_theorem}
Within the maximum waveform time span $\Delta T_{\max }^k$, the predictive solution ${\hat p_k}(x,t)$ of the boundary output waveform in each sub-domain is equal to the true solution $u(x,t)$ of the original domain $\omega$.
\end{theorem}

A simple proof for \cref{predict_theorem} is presented in \cref{section_appendix}.

\begin{theorem}[Waveform Time Span Theorem]
\label{velocity_theorem}
To assure the validation of RSWR, the maximum waveform time span should be less than $\Delta \tau/2$, 
where $\Delta \tau$ is defined as the minimum time that the wave costs to transmit through the overlapping region. 

\end{theorem}

Assume the width of the overlapping domain is  $\Delta {X_{overlap}}$, and the wave velocity of wave equation Eq. \cref{wave_equation} is $a$, then the waveform Maximum time span $\Delta {T_{\max }}$ satisfies:
\begin{equation}
\Delta {T_{\max }} < \frac{{\Delta {X_{overlap}}}}{2a}
\end{equation}


\subsection{Update}
\label{subsection_update}

This section updates the boundary input waveform and calculate the true solution for subdomains.

According to \cref{predict_theorem}, for Eq. \cref{sub_equation_predict_1}, the conclusion is as below:

\begin{equation}
\left\{
\begin{array}{l}
\frac{{\partial {p_1}(x,t)}}{{\partial x}}{|_{x = {X_{2A}}}} = \frac{{\partial {{\hat p}_1}(x,t)}}{{\partial x}}{|_{x = {X_{2A}}}} = \frac{{\partial u(x,t)}}{{\partial x}}{|_{x = {X_{2A}}}}, \\
\quad \quad \forall t \in [T_{start}^k,T_{start}^k + \Delta T_{\max }^k],k = 1 \\
\end{array}
\right.
\end{equation}

Similarly, for Eq. \cref{sub_equation_predict_2}, we have:

\begin{equation}
\left\{
\begin{array}{l}
\frac{{\partial {p_2}(x,t)}}{{\partial x}}{|_{x = {X_{1B}}}} = \frac{{\partial {{\hat p}_2}(x,t)}}{{\partial x}}{|_{x = {X_{1B}}}} = \frac{{\partial u(x,t)}}{{\partial x}}{|_{x = {X_{1B}}}}, \\
\quad \quad \forall t \in [T_{start}^k,T_{start}^k + \Delta T_{\max }^k],k = 1
\end{array}
\right.
\end{equation}

As the result, the boundary output waveform of Eq. \cref{sub_equation_1} is solved by Eq. \cref{sub_equation_predict_1}, and Eq. \cref{sub_equation_1} is updated as Eq. \cref{sub_equation_solve_1}:

\begin{equation}
\label{sub_equation_solve_1}
\left\{
\begin{array}{l}
 \frac{{{\partial ^2}{p_1}(x,t)}}{{\partial {t^2}}} - {a^2}\frac{{{\partial ^2}{p_1}(x,t)}}{{\partial {x^2}}} = 0,a > 0, x \in {\Omega _1},{\Omega _1} = \left[ {{X_{1A}},{X_{1B}}} \right],\\
 \quad \quad t \in [T_{start}^k,T_{start}^k + \Delta T_{\max }^k],k = 1 \\ 
 {p_1}(x,t){|_{t = T_{start}^k}} = u(x,t){|_{t = T_{start}^k}},x \in {\Omega _1} \\ 
 \frac{{\partial {p_1}(x,t)}}{{\partial t}}{|_{_{t = T_{start}^k}}} = \frac{{\partial u(x,t)}}{{\partial t}}{|_{_{t = T_{start}^k}}},x \in {\Omega _1} \\ 
 {p_1}(x,t){|_{x = {X_{1A}}}} = {f_A}(t),t \in [T_{start}^k,T_{start}^k + \Delta T_{\max }^k] \\ 
 \frac{{\partial {p_1}(x,t)}}{{\partial x}}{|_{x = {X_{1B}}}} = \frac{{\partial {{\hat p}_2}(x,t)}}{{\partial x}}{|_{x = {X_{1B}}}},t \in [T_{start}^k,T_{start}^k + \Delta T_{\max }^k]
 \end{array}
\right.
\end{equation}

Eq. \cref{sub_equation_solve_1} can be solved and obtain: \[{p_1}(x,t),x \in {\Omega _1},t \in [T_{start}^k,T_{start}^k + \Delta T_{\max }^k]\]

Similarly, Eq. \cref{sub_equation_predict_2} is updated as Eq. \cref{sub_equation_solve_2}:

\begin{equation}
\label{sub_equation_solve_2}
\left\{
\begin{array}{l}
 \frac{{{\partial ^2}{p_2}(x,t)}}{{\partial {t^2}}} - {a^2}\frac{{{\partial ^2}{p_2}(x,t)}}{{\partial {x^2}}} = 0,a > 0,x \in {\Omega _2},{\Omega _2} = \left[ {{X_{2A}},{X_{2B}}} \right],\\
\quad \quad t \in [T_{start}^k,T_{start}^k + \Delta T_{\max }^k],k = 1 \\ 
 {p_2}(x,t){|_{t = T_{start}^k}} = u(x,t){|_{t = T_{start}^k}},x \in {\Omega _2} \\ 
 \frac{{\partial {p_2}(x,t)}}{{\partial t}}{|_{t = T_{start}^k}} = \frac{{\partial u(x,t)}}{{\partial t}}{|_{t = T_{start}^k}},x \in {\Omega _2} \\ 
 \frac{{\partial {p_2}(x,t)}}{{\partial x}}{|_{x = {X_{2A}}}} = \frac{{\partial {{\hat p}_2}(x,t)}}{{\partial x}}{|_{x = {X_{2A}}}},t \in [T_{start}^k,T_{start}^k + \Delta T_{\max }^k] \\ 
 {p_2}(x,t){|_{x = {X_{2B}}}} = {f_B}(t),t \in [T_{start}^k,T_{start}^k + \Delta T_{\max }^k] \\ 
 \end{array}
\right.
\end{equation}

Eq. \cref{sub_equation_solve_2} can be solved and get: \[{p_2}(x,t),x \in {\Omega _2},t \in [T_{start}^k,T_{start}^k + \Delta T_{\max }^k],k = 1\]

Consequently, the true solution of Eq. \cref{wave_equation}:

\[u(x,t),  x \in \Omega,  t \in [T_{start}^k,T_{start}^k + \Delta T_{\max }^k],k = 1\] 

is the combination of the solution of Eq. \cref{sub_equation_solve_1} and Eq. \cref{sub_equation_solve_2}:

\begin{equation}
u(x,t) = \left\{ \begin{array}{l}
 {p_1}(x,t),x \in {\Omega _1} \\ 
 {p_2}(x,t),x \in {\Omega _2} \wedge x \notin {\Omega _1} \\ 
 \end{array} \right.,t \in [T_{start}^k,T_{start}^k + \Delta T_{\max }^k],k = 1
\end{equation}

\subsection{Redo the loop}
\label{subsection_redo}

Update the start time and initial condition, repeat Section \ref{subsection_split}, \ref{subsection_predict} and \ref{subsection_select}.

Set the new start time as:

\[T_{start}^k = T_{stop}^{k - 1} = T_{start}^{k - 1} + \Delta T_{\max }^{k - 1},k = 2\]
  
For Eq. \cref{sub_equation_predict_1}, set the initial condition at  $t = T_{start}^k, k = 2$ as:

\begin{equation}
\begin{array}{l}
 {p_1}(x,t){|_{t = T_{start}^k}} = u(x,t){|_{t = T_{start}^k}},x \in {\Omega _1}, k = 2 \\ 
 \frac{{\partial {p_1}(x,t)}}{{\partial t}}{|_{t = T_{start}^k}} = \frac{{\partial u(x,t)}}{{\partial t}}{|_{t = T_{start}^k}},x \in {\Omega _1}, k = 2 \\ 
 \end{array}
\end{equation}

Similarly, for Eq. \cref{sub_equation_predict_2}, set the initial condition at $t = T_{start}^k, k = 2$ as:

\begin{equation}
\begin{array}{l}
 {p_2}(x,t){|_{t = T_{start}^k}} = u(x,t){|_{t = T_{start}^k}}, x \in {\Omega _2}, k=2 \\ 
 \frac{{\partial {p_2}(x,t)}}{{\partial t}}{|_{t = T_{start}^k}} = \frac{{\partial u(x,t)}}{{\partial t}}{|_{t = T_{start}^k}},x \in {\Omega _2}, k=2 \\ 
 \end{array}
\end{equation}

Then, set the new predictive time span:

\[\Delta \hat T_{predict}^k = (1 + \beta )*\Delta \hat T_{\max }^{k - 1},k = 2, \beta  = 0.1\]

and redo Section \ref{subsection_predict} to find the new max time span $\Delta {T_{\max }}$, and redo Section \ref{subsection_select} to get the true solution of Eq. \cref{wave_equation} in the new time span:

\[u(x,t),x \in \Omega ,t \in [T_{start}^k,T_{start}^k + \Delta T_{\max }^k],k = 2\]

Repeat the above procedure by loop, the original wave equation Eq. \cref{wave_equation} is distributed solved by decomposing into 2 overlapped domains.

\section{Algorithm Extended for $N > 2$}
\label{section_algorithm_n}

Assume the original wave equation is split into $N$ overlapped subdomains by RSWR, by using $1$-Dimmension partition strategy.

\begin{definition}[Global Maximum Waveform Time Span]
  For each two domain  ${\Omega _i}$ and ${\Omega _j}$, there will be a maximum waveform time span $\Delta T_{\max }^k(i,j)$, which satisfies:
\end{definition}

\begin{equation} 
\left\{
\begin{array}{l}
\Delta T_{\max }^k(i,j) > 0, \quad if\;{\Omega _i} \cap {\Omega _j} \ne \emptyset  \\ 
\Delta T_{\max }^k(i,j) =  + \infty, \quad if\;{\Omega _i} \cap {\Omega _j} = \emptyset  \\ 
i = 1, \cdots ,N, \quad j = 1, \cdots ,N,i \ne j \\ 
k = 1, \cdots , + \infty  \\ 
\end{array} 
\right.
\end{equation}
 
Therefore, the global maximum waveform time span $\Delta T_{\max ,global}^k$ would be:
\begin{equation}
\Delta T_{\max ,global}^k = \min (\Delta T_{\max }^k(i,j)),i = 1, \cdots ,N,j = 1, \cdots ,N,i \ne j
\end{equation}

Consequently, \cref{predict_theorem} is extended into \cref{global_max_time_span_theorem}, 
where $p_i^k(x,t)$ is the true solution for sub-domain ${\Omega _i}$, 
and $\hat p_i^k(x,t)$ is the predictive solution for ${\Omega _i}$.

\begin{theorem}[Global Maximum Waveform Time Span Theorem]
\label{global_max_time_span_theorem}
Within the global maximum waveform time span $\Delta T_{\max ,global}^k$, the predictive solution of the boundary output waveform in each sub-domain ${\Omega _i}$ is equal to the true solution of ${\Omega _i}$.
\end{theorem}

According to \cref{global_max_time_span_theorem}: 

\begin{equation}
\left\{
\begin{array}{l}
\frac{{\partial p_i^k(x,t)}}{{\partial x}}{|_{x = {X_{OB,i}}}} = 
\frac{{\partial \hat p_i^k(x,t)}}{{\partial x}}{|_{x = {X_{OB,i}}}} = 
\frac{{\partial u(x,t)}}{{\partial x}}{|_{x = {X_{OB,i}}}}, \\
\quad \quad x \in {\Omega _i}, i = 1,\cdots,N - 1, \\
\quad \quad t \in [T_{start}^k,T_{start}^k + \Delta T_{\max ,global}^k] \\
\end{array}
\right.
\end{equation}

where $x = {X_{OB,i}}$ is the output boundary in sub-domain ${\Omega _i}$.

Based on \cref{global_max_time_span_theorem}, the predict-select-update strategy of RSWR is valid for $N>2$.

\section{Experiments}
\label{section_experiment}

\subsection{$N=2$} 
\label{subsection_experiment_1}

\begin{figure}[htb]
  \centering
  \label{fig_n_2_true_solution}
  \includegraphics[width=\textwidth]{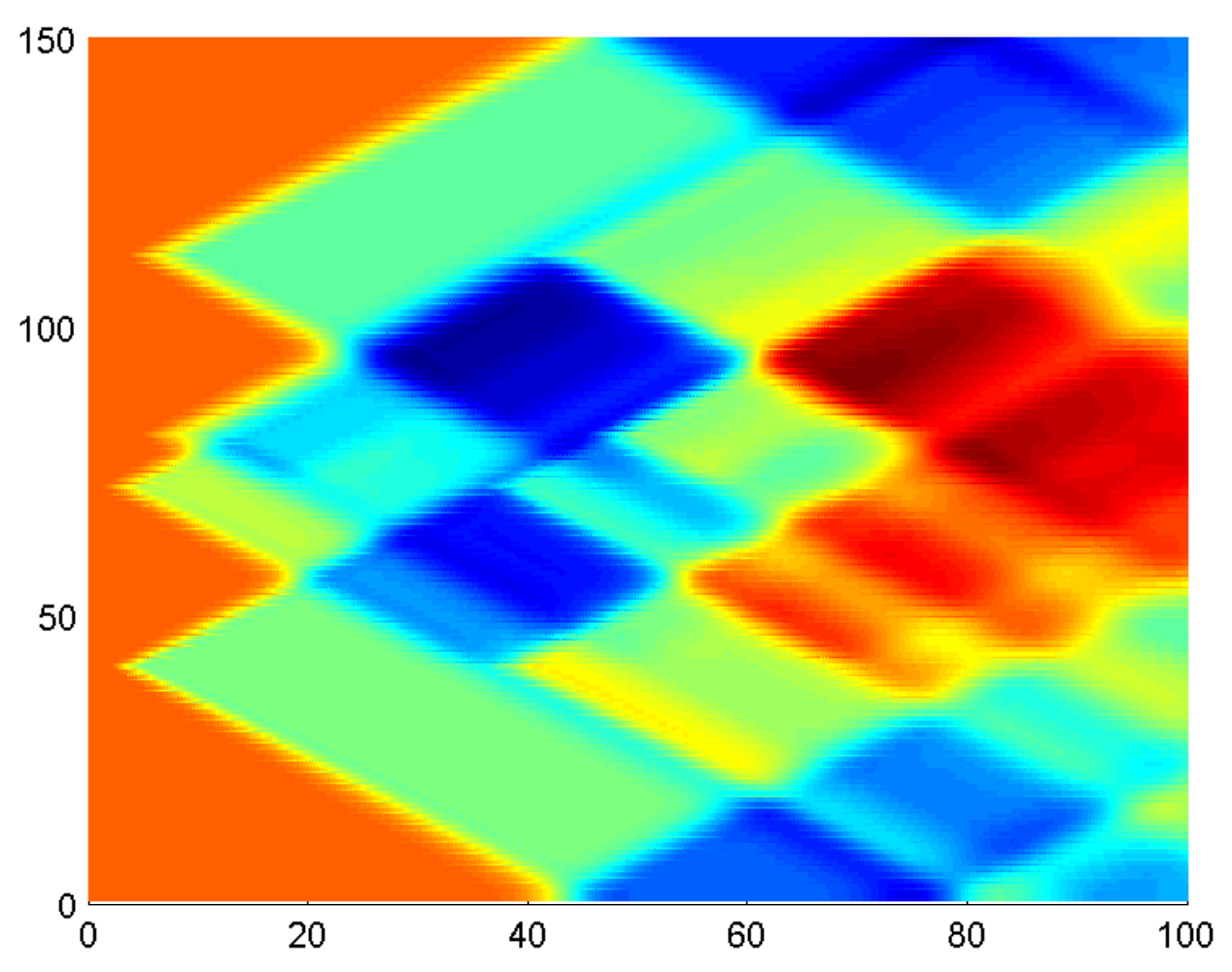}
  \caption{Numerical true solution of Eq. \cref{wave_equation}, $N=1$.}
\end{figure}

\begin{figure}[htb]
  \centering
  \label{fig_n_2_RS_solution}
  \includegraphics[width=\textwidth]{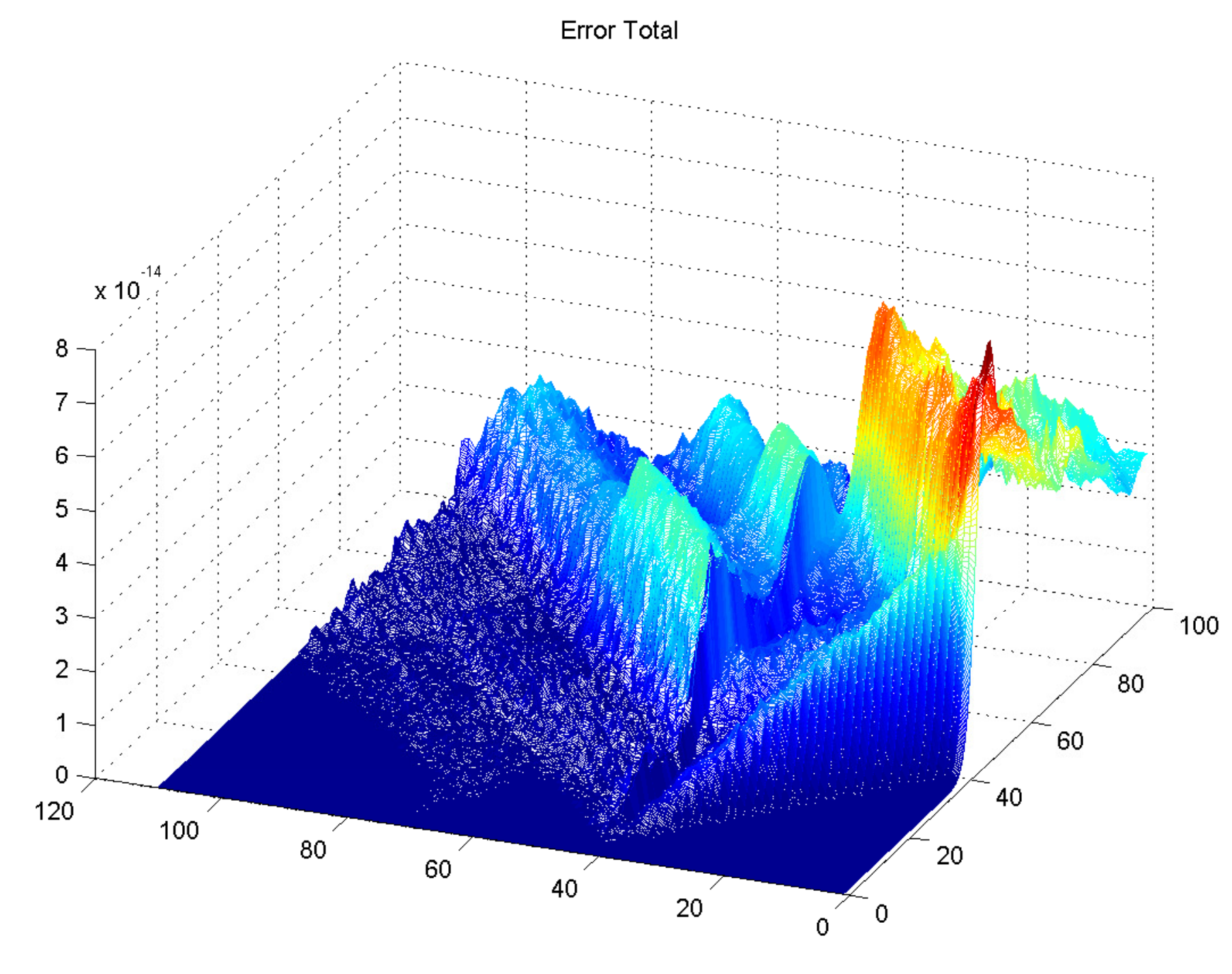}
  \caption{Distributed numerical solution of Eq. \cref{wave_equation} by RSWR, $N=2$.}
\end{figure}

The 1-D wave equation \cref{wave_equation} is inserted with $N=2$ pulse sources, the numerical result for Eq. \cref{wave_equation} is shown as \cref{fig_n_2_true_solution}.
Then we split Eq. \cref{wave_equation} into $N=2$ sub-domains by 1-D partitioning, and use RSWR to calculate them distributedly. The error of RSWR is shown as \cref{fig_n_2_RS_solution}.

\subsection{$N=10$}
\label{subsection_experiment_2}
\begin{figure}[htb]
  \centering
  \label{fig_n_10_true_solution}
  \includegraphics[width=\textwidth]{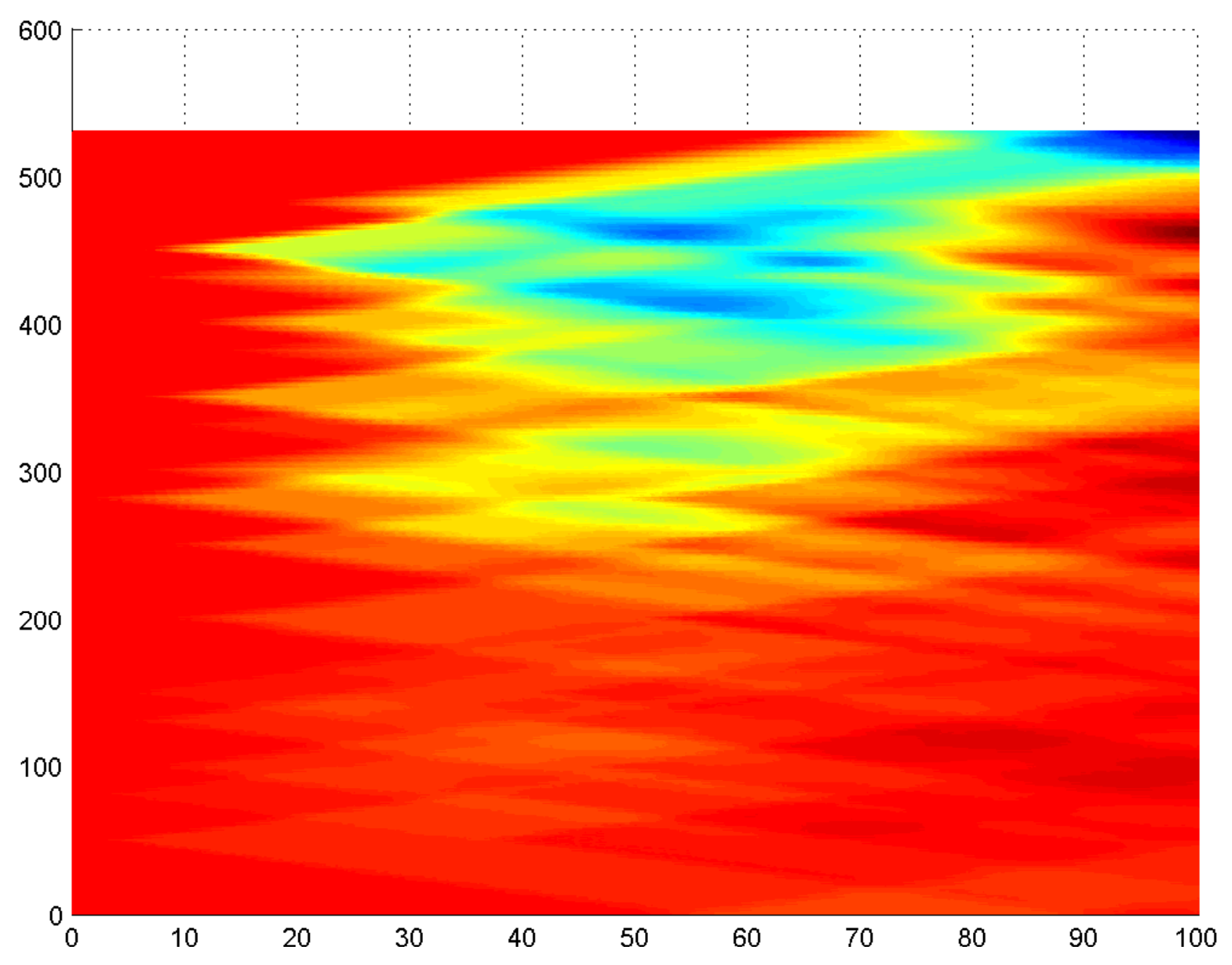}
  \caption{Numerical true solution of Eq. \cref{wave_equation}, $N=1$.}
\end{figure}

\begin{figure}[htb]
  \centering
  \label{fig_n_10_RS_solution}
  \includegraphics[width=\textwidth]{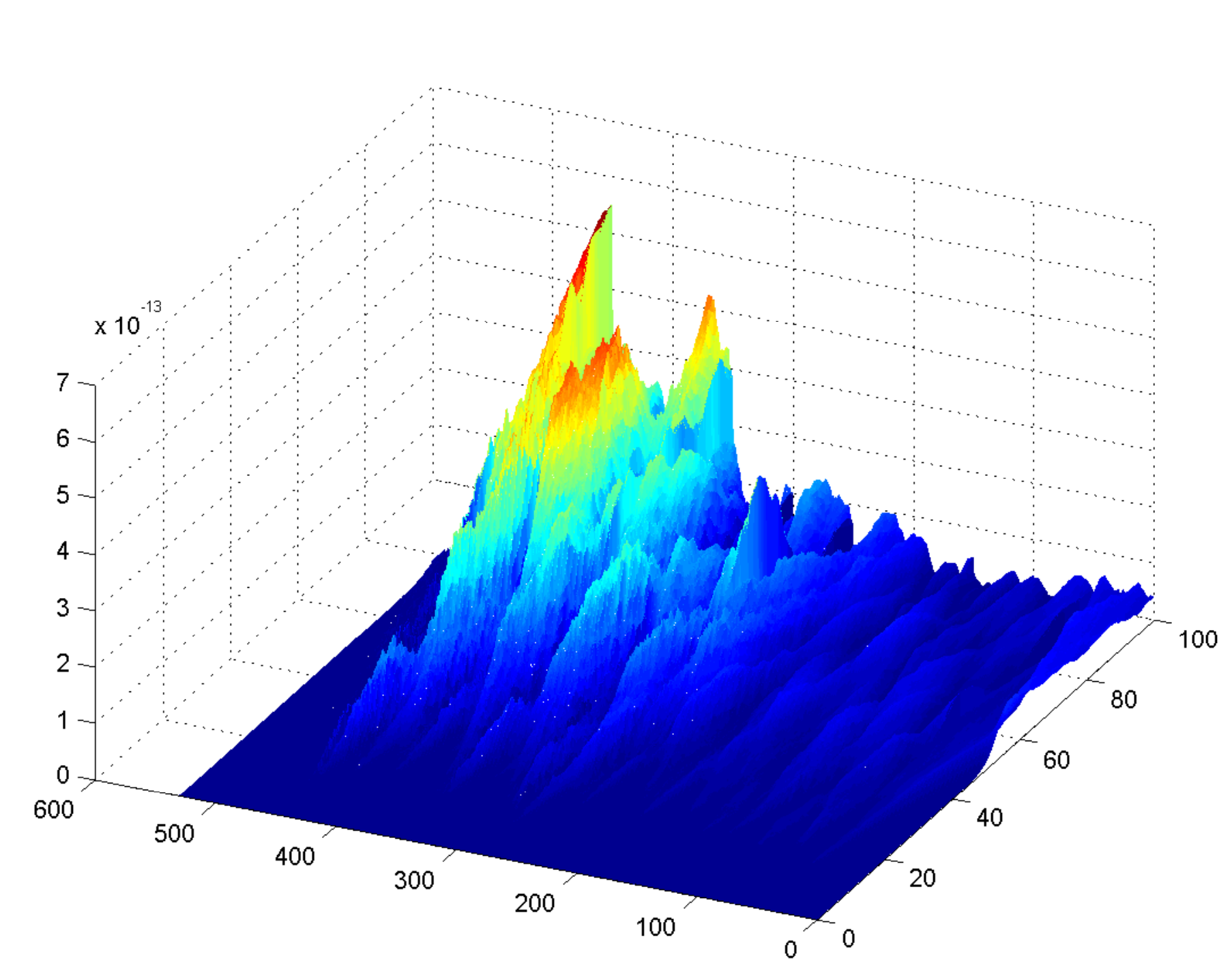}
  \caption{Distributed numerical solution of Eq. \cref{wave_equation} by RSWR, $N=10$.}
\end{figure}

The 1-D wave equation \cref{wave_equation} is inserted with $N=10$ pulse sources, the numerical result for Eq. \cref{wave_equation} is shown as \cref{fig_n_10_true_solution}.
Then we split Eq. \cref{wave_equation} into $N=10$ sub-domains by 1-D partitioning, and use RSWR to calculate them distributedly. The error of RSWR is shown in \cref{fig_n_10_RS_solution}.

\section{Conclusions}
\label{section_conclusion}
This paper proposes Relative Schwarz Waveform Relaxation algorithm (RSWR), 
which is a non-iterative overlapping SWR for wave equation. 
RSWR is able to achieve high accuracy result by using the predict-select-update strategy,
and it does not need preconditioner.
Experiments show that the accuracy of RSWR is good, 
and RSWR is potential to be scalable. 
Since RSWR is a non-iterative algorithm, it is potential to be fast.
Therefore, it would be meaningful to implement RSWR on supercomputers 
to solve large physical problem. 

Beyond wave equation, RSWR is potential to solve linear and nonlinear hyperbolic partial differential equation distributedly, which requires further study.

\appendix
\section{Proof for \cref{predict_theorem}}
\label{section_appendix} 
To be complemented.

\section*{Acknowledgments}
Fei Wei would like to thank his parents, Chengluan Wei and Qing Li, for their support and comfort during his hard times. Fei Wei would also like to thank his wife, Anna Zhao, who encourages him to continue this long term research which started at 2006. Further, Fei Wei would like to thank his friends, Yao Yu, Xiaoyang Yang, Qi Wei, Peng Zhang, for their help.

\bibliography{main}
\end{document}